\documentclass[12pt]{article}
\usepackage{amsfonts}
\usepackage{myart}

\normalbaselines
\input tcilatex

\begin{document}

\title{Euclidean geometry as algorithm for construction of generalized
geometries.}
\author{Yuri A. Rylov}
\date{}
\maketitle

\begin{abstract}
It is shown that the generalized geometries may be obtained as a deformation
of the proper Euclidean geometry. Algorithm of construction of any
proposition $\mathcal{S}$ of the proper Euclidean geometry $\mathcal{E}$ may
be described in terms of the Euclidean world function $\sigma _{\mathrm{E}}$
in the form $\mathcal{S}\left( \sigma _{\mathrm{E}}\right) $. Replacing the
Euclidean world function $\sigma _{\mathrm{E}}$ by the world function $%
\sigma $ of the geometry $\mathcal{G}$, one obtains the corresponding
proposition $\mathcal{S}\left( \sigma \right) $ of the generalized geometry $%
\mathcal{G}$. Such a construction of the generalized geometries (known as
T-geometries) uses well known algorithms of the proper Euclidean geometry
and nothing besides. This method of the geometry construction is very simple
and effective. Using T-geometry as the space-time geometry, one can
construct the deterministic space-time geometries with primordially
stochastic motion of free particles and geometrized particle mass. Such a
space-time geometry defined properly (with quantum constant as an attribute
of geometry) allows one to explain quantum effects as a result of the
statistical description of the stochastic particle motion (without a use of
quantum principles).
\end{abstract}

\section{Introduction}

The proper Euclidean geometry has been constructed by Euclid many years ago.
The Euclidean geometry may be considered as a set of many algorithms, which
are necessary for construction of geometrical objects $\mathcal{O}_{\mathrm{E%
}}$ and relations $\mathcal{R}_{\mathrm{E}}$ between them. These algorithms
of the geometrical objects construction were obtained by means of logical
reasonings from fundamental propositions (axioms) of Euclidean geometry. Any
algorithm $\mathcal{A}_{\mathrm{E}}$ of construction of some geometrical
object $\mathcal{O}_{\mathrm{E}}$ may be considered as an operator $O$,
acting on the set $\Omega $ of points $P$, where the geometry is
constructed. The objects, which undergo the action of operators, are called
operands. In the given case the points $P$ of the set $\Omega $ are operands
of algorithms of the Euclidean geometry construction.

There is a necessity in construction of generalized (non-Euclidean)
geometries, which are distinguished from the proper Euclidean geometry in
different aspects. For instance, a generalized geometry is necessary for
description of properties of the real space-time.

In the beginning of 20th century there existed the problem of motion of
bodies with very high speed. This problem was solved by a construction of
the relativity theory, where the direct product of time and the Euclidean
space was substituted by the Minkowski geometry of the space-time. Another
problem of the beginning of 20th century: stochastic motion of particle with
small mass also should be solved by means of a modification of the
space-time geometry. The fact is that the motion of free particles is
determined \textit{only by the properties of the space-time} (or by
properties of its geometry). Intensity of the stochastic component of the
particle motion depends on the particle mass. This intensity is very small
for particles of large mass, and it is essential for particles of small
mass. It means that in the corresponding space-time geometry the free
particle motion is to be primordially stochastic, and the particle mass is
to be geometrized, to take into account the influence of the mass upon the
free motion. In the beginning of 20th century there was not such a geometry
in the framework of the Riemannian geometries. Furthermore, one could not
imagine such a deterministic space-time geometry, where the motion of free
particles be stochastic. One tried to solve the problem, using stochastic
space-time geometry \cite{M51,B70}. However, the stochastic geometry is not
a geometry in the precise value of the word. The stochastic geometry is the
usual space-time geometry (for instance, Minkowski geometry) with some
additional structure, given on the space-time.

As a result the problem of the stochastic particle motion had not been
solved on the fundamental level (by a modification of the space-time
geometry). Instead, the stochastic manifestation of the space-time geometry
properties was prescribed to the free particle in the form of quantum
properties, generated by the quantum principles. Introduction of quantum
principles admitted one to explain all nonrelativistic quantum effects.
However, the relativistic quantum effects cannot be described on the basis
of the quantum principles, because quantum principles are nonrelativistic.
To expand quantum theory to relativistic phenomena, one needs to return on
the fundamental level, i. e. to the true space-time geometry, which is valid
for both relativistic and nonrelativistic phenomena.

Conventional method of the generalized geometry construction is a
modification of the Euclidean algorithms. One changes some Euclidean axioms.
Besides, one uses deformation of the Euclidean space. The infinitesimal
Euclidean distance $dS_{\mathrm{E}}=\sqrt{g_{\mathrm{E}ik}dx^{i}dx^{k}}$ is
substituted by means of the Riemannian one $dS_{\mathrm{R}}=\sqrt{g_{\mathrm{%
R}ik}dx^{i}dx^{k}}$. In general, the Euclidean algorithm of the geometry
construction is replaced by another one. The new algorithm is used for
construction of the generalized geometry. To construct the generalized
geometry (for instance, the Riemannian geometry), one repeats many
(sometimes almost all) Euclidean logical reasonings . Such a method is very
complicated. Besides, it is necessary, that independent modifications of
different Euclidean axioms be compatible between themselves.

There is another method of the Euclidean geometry generalization. We do not
change the Euclidean algorithms of the geometry construction. However, we
replace operands of the Euclidean algorithms. We use the world function $%
\sigma $ as the operand of the Euclidean algorithms instead of the set $%
\Omega $ of points $P$, which are usually is used as operands of the
Euclidean algorithms.

The world function $\sigma \left( P,Q\right) =\frac{1}{2}\rho ^{2}\left(
P,Q\right) $, where $\rho \left( P,Q\right) $ is the distance between the
points $P$,$Q\in \Omega $. (We do not use the term metric for $\rho $,
because the metric is supposed to be restricted by constraint of positivity
and by the triangle axiom, whereas the distance $\rho $ is free of these
constraints.) The world function $\sigma $ defined by the relation 
\begin{equation}
\sigma :\qquad \Omega \otimes \Omega \rightarrow \mathbb{R},\qquad \sigma
\left( P,P\right) =0,\qquad \forall P\in \Omega  \label{d1.1}
\end{equation}
The world function was introduced by Synge \cite{S60}, and it is a very
important quantity in the proper Euclidean geometry. But the world function $%
\sigma $ is not only an important quantity. It is the \textit{only important}
quantity in the proper Euclidean geometry. It means that the proper
Euclidean geometry and all its algorithms may be described in terms and only
in terms of the Euclidean world function $\sigma _{\mathrm{E}}$. The
property of a geometry to be described completely in terms and only in terms
of the world function is called the $\sigma $-immanence property.
Description of a geometry in terms of the world function is the $\sigma $%
-immanent description of the geometry.

The $\sigma $-immanence property of the Euclidean geometry was discovered
rather recently \cite{R90,R01}. It has been proved that the Euclidean
geometry can be presented in terms and only in terms of the function $\sigma
_{\mathrm{E}}$, provided the function $\sigma _{\mathrm{E}}$ satisfies a
series of constraints, written in terms of $\sigma _{\mathrm{E}}$. By
definition, \textit{any geometry is a totality of all geometric objects }$%
\mathcal{O}$\textit{\ and of all relations }$\mathcal{R}$\textit{\ between
them}. The $\sigma $-immanence of the proper Euclidean geometry means that
any geometric object $\mathcal{O}_{\mathrm{E}}$ and any relation $\mathcal{R}%
_{\mathrm{E}}$ of the Euclidean geometry $\mathcal{G}_{\mathrm{E}}$ can be
presented in terms of the Euclidean world function $\sigma _{\mathrm{E}}$ in
the form $\mathcal{O}_{\mathrm{E}}\left( \sigma _{\mathrm{E}}\right) $ and $%
\mathcal{R}_{\mathrm{E}}\left( \sigma _{\mathrm{E}}\right) $.

Let us suppose that a geometry $\mathcal{G}$ has the property of the $\sigma 
$-immanence. Then the geometry $\mathcal{G}$ may be constructed as a result
of a deformation of the proper Euclidean geometry $\mathcal{G}_{\mathrm{E}}$%
. Indeed, the proper Euclidean geometry $\mathcal{G}_{\mathrm{E}}$ is the
totality of geometrical objects $\mathcal{O}_{\mathrm{E}}\left( \sigma _{%
\mathrm{E}}\right) $ and relations $\mathcal{R}_{\mathrm{E}}\left( \sigma _{%
\mathrm{E}}\right) $. We produce the change

\begin{equation}
\sigma _{\mathrm{E}}\rightarrow \sigma ,\qquad \mathcal{O}_{\mathrm{E}%
}\left( \sigma _{\mathrm{E}}\right) \rightarrow \mathcal{O}_{\mathrm{E}%
}\left( \sigma \right) ,\qquad \mathcal{R}_{\mathrm{E}}\left( \sigma _{%
\mathrm{E}}\right) \rightarrow \mathcal{R}_{\mathrm{E}}\left( \sigma \right)
\label{h1.2}
\end{equation}
Then totality of geometrical objects $\mathcal{O}_{\mathrm{E}}\left( \sigma
\right) $, relations $\mathcal{R}_{\mathrm{E}}\left( \sigma \right) $ and
the world function $\sigma $ form the generalized geometry $\mathcal{G}$.
Any generalized geometry, obtained by the method (\ref{h1.2}) will be
referred to as a tubular geometry (T-geometry) \cite{R90,R01,R002}.

Note that the geometrical objects $\mathcal{O}_{\mathrm{E}}\left( \sigma _{%
\mathrm{E}}\right) $, $\mathcal{O}_{\mathrm{E}}\left( \sigma \right) $ and
relations $\mathcal{R}_{\mathrm{E}}\left( \sigma _{\mathrm{E}}\right) $, $%
\mathcal{R}_{\mathrm{E}}\left( \sigma \right) $ are constructed by means the
same algorithm of the Euclidean geometry construction. Only operands $\sigma
_{\mathrm{E}}$ and $\sigma $ are different in the two cases. Thus, different
generalized geometries $\mathcal{G}_{1}$ and $\mathcal{G}_{2}$ are
distinguished by the value of their operands $\sigma _{1}$ and $\sigma _{2}$%
, but not by the algorithms of their construction. It is very convenient,
because a change of operand is more simpler, than a change of an algorithm.
We may not care for compatibility of different changes of the Euclidean
algorithms, as it takes place at the conventional approach to a construction
of a generalized geometry, when we change the algorithm at the constant
operand. Besides, we can evaluate the power of the set of generalized
geometries, which can be considered as a power of the set of functions of
two arguments. The power of this set is much more, than the power of the set
of all Riemannian geometries, which can be evaluated as the power of the set
of several functions of one argument. For instance, the set of all
homogeneous isotropic T-geometries is determined by a function of one
argument, whereas the set of all homogeneous isotropic Riemannian geometries
is labelled by the dimension $n$ and the index $\nu =1,2,...n$.

From practical point of view, the T-geometries are interesting in the
relation, that they contains such space-time geometries, where the motion of
free particle is primordially stochastic. Furthermore, the space-time
T-geometry with stochastic motion of free particle is a general case,
whereas the space-time T-geometry with deterministic motion of free particle
is a degenerate case.

Further we shall formulate the principal theorem on the $\sigma $-immanence
of the proper Euclidean geometry. But at first, we put the cases of $\sigma $%
-immanent description of objects of the proper Euclidean geometry

Let $\mathcal{R}_{\mathrm{E}}\left( \sigma _{\mathrm{E}}\right) $ be the
scalar product $\left( \mathbf{P}_{0}\mathbf{P}_{1}.\mathbf{P}_{0}\mathbf{P}%
_{2}\right) _{\mathrm{E}}$ of two vectors $\mathbf{P}_{0}\mathbf{P}_{1}$, $%
\mathbf{P}_{0}\mathbf{P}_{2}$ in $\mathcal{G}_{\mathrm{E}}$. It can be
written in the $\sigma $-immanent form (i.e. in terms of the world function $%
\sigma _{\mathrm{E}}$) 
\begin{equation}
\mathcal{R}_{\mathrm{E}}\left( \sigma _{\mathrm{E}}\right) :\qquad \left( 
\mathbf{P}_{0}\mathbf{P}_{1}.\mathbf{P}_{0}\mathbf{P}_{2}\right) _{\mathrm{E}%
}=\sigma _{\mathrm{E}}\left( P_{0},P_{1}\right) +\sigma _{\mathrm{E}}\left(
P_{0},P_{2}\right) -\sigma _{\mathrm{E}}\left( P_{1},P_{2}\right)
\label{h1.3}
\end{equation}
where index 'E' shows that the quantity relates to the Euclidean geometry.
It is easy to see that (\ref{h1.3}) is a corollary of the Euclidean
relations 
\begin{equation}
\left\vert \mathbf{P}_{0}\mathbf{P}_{1}\right\vert _{\mathrm{E}}^{2}=2\sigma
_{\mathrm{E}}\left( P_{0},P_{1}\right)  \label{h1.7}
\end{equation}
\begin{equation}
\left\vert \mathbf{P}_{1}\mathbf{P}_{2}\right\vert _{\mathrm{E}%
}^{2}=\left\vert \mathbf{P}_{0}\mathbf{P}_{2}-\mathbf{P}_{0}\mathbf{P}%
_{1}\right\vert _{\mathrm{E}}^{2}=\left\vert \mathbf{P}_{0}\mathbf{P}%
_{1}\right\vert _{\mathrm{E}}^{2}+\left\vert \mathbf{P}_{0}\mathbf{P}%
_{2}\right\vert _{\mathrm{E}}^{2}-2\left( \mathbf{P}_{0}\mathbf{P}_{1}.%
\mathbf{P}_{0}\mathbf{P}_{2}\right) _{\mathrm{E}}  \label{h1.8}
\end{equation}

According to (\ref{h1.2}) in the generalized geometry $\mathcal{G}$ we
obtain instead of (\ref{h1.3}) 
\begin{equation}
\mathcal{R}_{\mathrm{E}}\left( \sigma \right) :\qquad \left( \mathbf{P}_{0}%
\mathbf{P}_{1}.\mathbf{P}_{0}\mathbf{P}_{2}\right) =\sigma \left(
P_{0},P_{1}\right) +\sigma \left( P_{0},P_{2}\right) -\sigma \left(
P_{1},P_{2}\right)  \label{h1.4}
\end{equation}

Such a way of the generalized geometry construction is very simple. It does
not use any logical reasonings. It is founded on the supposition that 
\textit{any generalized geometry has the }$\sigma $\textit{-immanence
property}. It uses essentially the fact that the algorithms of the proper
Euclidean geometry construction are already known, and all necessary logical
reasonings has been already produced in the proper Euclidean geometry.

The application of the replacement (\ref{h1.2}) to the construction of a
generalized geometry will be referred to as the deformation principle. Any
change of distance $\rho $, or the world function $\sigma $ between the
points of the space $\Omega $ means a deformation of this space. We construe
the concept of deformation in a broad sense. The deformation may transform a
point into a surface and a surface into a point. The deformation may remove
some points of the Euclidean space, violating its continuity, or decreasing
its dimension. The deformation may add supplemental points to the space,
increasing its dimension. We may interpret any $\sigma $-immanent
generalization of the Euclidean geometry as its deformation. In other words,
the deformation principle is a very general method of the generalized
geometry construction.

Construction of a nonhomogeneous geometry on the axiomatic basis is
impossible practically, because there is a lot of different nonhomogeneous
geometries. It is very difficult to invent axiomatics for a nonhomogeneous
geometry, where identical objects have different properties in various
places. Besides, one cannot invent axiomatics for each of these geometries.
Thus, in reality there is no alternative to application of the deformation
principle at the construction of the nonhomogeneous generalized geometry.
The real problem consists in the sequential application of the deformation
principle. As far as the \textit{deformation principle alone is sufficient
for the construction of the T-geometry}, one may not use additional means of
the geometry construction. At the T-geometry construction we do not use
coordinate system and other means of description.

Mathematicians provide physicists with their geometrical construction, and
physicists believe that the space-time is continuous. Continuity of the
space-time cannot be tested experimentally, and the only reason of the
space-time continuity is the fact that mathematicians are able to construct
only continuous geometries, whereas they fail to construct discrete
geometries.

We have the same situation with the space dimension. Geometers consider the
dimension to be an inherent property of any geometry. They can imagine the $%
n $-dimensional Riemannian geometry, but they cannot imagine a geometry
without a dimension, or a geometry of an indefinite dimension. The reason of
these belief is the fact that the dimension of the manifold and its
continuity are the starting points of the Riemannian geometry construction,
and at this point one cannot separate the properties of the geometry from
the properties of the manifold.

In the T-geometry we deal only with the geometry in itself, because it does
not use any means of the description. As a result the T-geometry is
insensitive to continuity or discreteness of the space, as well as to its
dimension. Application of additional means of description can lead to
inconsistency and to a restriction of the list of possible T-geometries.

Any generalization of the proper Euclidean geometry is founded on some
property of the Euclidean geometry (or its objects). This property is
conserved in all generalized geometries, whereas other properties of the
Euclidean geometry are varied. Character and properties of the obtained
generalized geometry depend essentially on the choice of the conserved
property of the basic Euclidean geometry. For instance, the Riemannian
geometry is such a generalization of the Euclidean one, where the
one-dimensionality and continuity of the Euclidean straight line are
conserved, whereas its curvature and torsion are varied. The straight line
is considered to be the principal geometric object of the Euclidean
geometry, and one supposes that such properties of the Euclidean straight
line as continuity and one-dimensionality (absence of thickness) are to be
conserved at the generalization. It means that the continuity and
one-dimensionality of the straight line are to be the principal concepts of
the generalized geometry (the Riemannian geometry). In accordance with such
a choice of the conserved geometrical object one introduces the concept of
the curve $\mathcal{L}$ as a continuous mapping of a segment of the real
axis onto the space $\Omega $%
\begin{equation}
\mathcal{L}:\qquad \left[ 0,1\right] \rightarrow \Omega  \label{g1.1}
\end{equation}
To introduce the concept of the continuity, which is a basic concept of the
generalization, one introduces the topological space, the dimension of the
space $\Omega $ and other basic concepts of the Riemannian geometry, which
are necessary for construction of the Riemannian generalization of the
Euclidean geometry. However, these properties are not necessary for the
T-geometry construction.

The $\sigma $-immanence of the Euclidean geometry is a property of the whole
Euclidean geometry. Using the $\sigma $-immanence for generalization, we do
not impose any constraints on the single geometric objects of the Euclidean
geometry. As a result the $\sigma $-immanent generalization appears to be a
very powerful generalization. Besides, from the common point of view the
application of the whole geometry property for the generalization seems to
be more reasonable, than a use of the properties of a single geometric
object. Thus, using the property of the whole Euclidean geometry, the $%
\sigma $-immanent generalization seems to be more reasonable, than the
Riemannian generalization, which uses the properties of the Euclidean
straight line.

Now we list the most attractive features of the $\sigma $-immanent
generalization of the Euclidean geometry:

\begin{enumerate}
\item The $\sigma $-immanent generalization uses the $\sigma $-immanence,
which is a property of the Euclidean geometry as a whole (but not a property
of a single geometric object as it takes place at the Riemannian
generalization).

\item The $\sigma $-immanent generalization does not use any logical
construction, and the $\sigma $-immanent generalization is automatically as
consistent, as the Euclidean geometry, whose axiomatics is used implicitly.
In particular, the T-geometry does not contain any theorems. As a result the
main problem of the T-geometry is a correct $\sigma $-immanent description
of geometrical objects and relations of the Euclidean geometry. There are
some subtleties in such a $\sigma $-immanent description, which are
discussed below.

\item The $\sigma $-immanent generalization is a very powerful
generalization. It varies practically all properties of the Euclidean
geometry, including such ones as the continuity and the parallelism
transitivity, which are conserved at the Riemannian generalization.

\item The $\sigma $-immanent generalization allows one to use the
coordinateless description and to ignore the problems, connected with the
coordinate transformations as well as with the transformation of other means
of description.

\item The T-geometry may be used as the space-time geometry. In this case
the tubular character of straights explains freely the stochastic world
lines of quantum microparticles. Considering the quantum constant $\hbar $
as an attribute of the space-time geometry, one can obtain the quantum
description as the statistical description of the stochastic world lines 
\cite{R91}. Such a space-time geometry cannot be obtained in the framework
of the Riemannian generalization of the Euclidean geometry.

\item Practically all above mentioned properties of the $\sigma $-immanent
generalization are corollaries of the fact, that the world function $\sigma $
is an operand in the algorithms of the proper Euclidean geometry
construction.
\end{enumerate}

\section{Euclidean geometry in the $\protect\sigma $-immanent form}

\begin{definition}
The $\sigma $-space $V=\left\{ \sigma ,\Omega \right\} $ is the set $\Omega $
of points $P$ with the given world function $\sigma $
\begin{equation}
\sigma :\qquad \Omega \times \Omega \rightarrow \mathbb{R},\qquad \sigma
\left( P,P\right) =0,\qquad \forall P\in \Omega   \label{h1.1}
\end{equation}
\end{definition}

Let the proper Euclidean geometry be given on the set $\Omega $, and the
quantity 
\begin{equation}
\rho \left( P_{0},P_{1}\right) =\sqrt{2\sigma \left( P_{0},P_{1}\right) }%
,\qquad P_{0},P_{1}\in \Omega  \label{g2.1}
\end{equation}
be the Euclidean distance between the points $P_{0},P_{1}$.

Let the vector $\mathbf{P}_{0}\mathbf{P}_{1}\mathbf{=}\left\{
P_{0},P_{1}\right\} $ be the ordered set of two points $P_{0}$, $P_{1}$. The
point $P_{0}$ is the origin of the vector $\mathbf{P}_{0}\mathbf{P}_{1}$,
and the point $P_{1}$ is its end. The length $\left\vert \mathbf{P}_{0}%
\mathbf{P}_{1}\right\vert $ of the vector $\mathbf{P}_{0}\mathbf{P}_{1}$ is
defined by the relation 
\begin{equation}
\left\vert \mathbf{P}_{0}\mathbf{P}_{1}\right\vert ^{2}=2\sigma \left(
P_{0},P_{1}\right)  \label{g2.2}
\end{equation}
In the Euclidean geometry the scalar product $\left( \mathbf{P}_{0}\mathbf{P}%
_{1}.\mathbf{P}_{0}\mathbf{P}_{2}\right) $ of two vectors $\mathbf{P}_{0}%
\mathbf{P}_{1}$ and $\mathbf{P}_{0}\mathbf{P}_{2}$, having the common origin 
$P_{0}$, is expressed by the relation (\ref{h1.4}) 
\begin{equation}
\left( \mathbf{P}_{0}\mathbf{P}_{1}.\mathbf{P}_{0}\mathbf{P}_{2}\right)
=\sigma \left( P_{0},P_{1}\right) +\sigma \left( P_{1},P_{0}\right) -\sigma
\left( P_{1},P_{2}\right)  \label{g2.3}
\end{equation}

It follows from the expression (\ref{g2.3}), written for scalar products $%
\left( \mathbf{P}_{0}\mathbf{P}_{1}.\mathbf{P}_{0}\mathbf{Q}_{1}\right) $
and $\left( \mathbf{P}_{0}\mathbf{P}_{1}.\mathbf{P}_{0}\mathbf{Q}_{0}\right) 
$, and from the properties of the scalar product in the Euclidean space,
that the scalar product $\left( \mathbf{P}_{0}\mathbf{P}_{1}.\mathbf{Q}_{0}%
\mathbf{Q}_{1}\right) $ of two vectors $\mathbf{P}_{0}\mathbf{P}_{1}$ and $%
\mathbf{Q}_{0}\mathbf{Q}_{1}$ can be written in the $\sigma $-immanent form 
\begin{eqnarray}
\left( \mathbf{P}_{0}\mathbf{P}_{1}.\mathbf{Q}_{0}\mathbf{Q}_{1}\right)
&=&\left( \mathbf{P}_{0}\mathbf{P}_{1}.\mathbf{P}_{0}\mathbf{Q}_{1}\right)
-\left( \mathbf{P}_{0}\mathbf{P}_{1}.\mathbf{P}_{0}\mathbf{Q}_{0}\right) = 
\nonumber \\
&&\sigma \left( P_{0},Q_{1}\right) +\sigma \left( P_{1},Q_{0}\right) -\sigma
\left( P_{0},Q_{0}\right) -\sigma \left( P_{1},Q_{1}\right)  \label{b1.1a}
\end{eqnarray}

Let $\mathbf{P}_{0}\mathbf{P}_{1}$, $\mathbf{P}_{0}\mathbf{P}_{2}$,...$%
\mathbf{P}_{0}\mathbf{P}_{n}$ be $n$ vectors in the Euclidean space. The
necessary and sufficient condition of their linear dependence is 
\begin{equation}
F_{n}\left( \mathcal{P}^{n}\right) \equiv \det \left\vert \left\vert \left( 
\mathbf{P}_{0}\mathbf{P}_{i}.\mathbf{P}_{0}\mathbf{P}_{k}\right) \right\vert
\right\vert =0,\qquad i,k=1,2,..n,\qquad \mathcal{P}^{n}=\left\{
P_{0},P_{1},...P_{n}\right\}  \label{g2.4}
\end{equation}
where $F_{n}\left( \mathcal{P}^{n}\right) \equiv \det \left\vert \left\vert
\left( \mathbf{P}_{0}\mathbf{P}_{i}.\mathbf{P}_{0}\mathbf{P}_{k}\right)
\right\vert \right\vert $ is the Gram's determinant, constructed of the
scalar products of vectors.

Let us formulate the theorem on the $\sigma $-immanence of the Euclidean
geometry.

\begin{theorem}
The $\sigma $-space $V=\left\{ \sigma ,\Omega \right\} $ is the $n$%
-dimensional proper Euclidean space, if and only if the world function $%
\sigma $ satisfies the following conditions, written in terms of the world
function $\sigma $.
\end{theorem}

\textit{I. Condition of symmetry: } 
\begin{equation}
\sigma \left( P,Q\right) =\sigma \left( Q,P\right) ,\qquad \forall P,Q\in
\Omega  \label{a1.4}
\end{equation}

\textit{II. Definition of the dimension: } 
\begin{equation}
\exists \mathcal{P}^{n}\equiv \left\{ P_{0},P_{1},...P_{n}\right\} \subset
\Omega ,\qquad F_{n}\left( \mathcal{P}^{n}\right) \neq 0,\qquad F_{k}\left( {%
\Omega }^{k+1}\right) =0,\qquad k>n  \label{g2.5}
\end{equation}
\textit{where }$F_{n}\left( \mathcal{P}^{n}\right) $\textit{\ is the Gram's
determinant (\ref{g2.4}). Vectors }$\mathbf{P}_{0}\mathbf{P}_{i}$\textit{, }$%
\;i=1,2,...n$\textit{\ are basic vectors of the rectilinear coordinate
system }$K_{n}$\textit{\ with the origin at the point }$P_{0}$\textit{, and
the metric tensors }$g_{ik}\left( \mathcal{P}^{n}\right) $\textit{, }$%
g^{ik}\left( \mathcal{P}^{n}\right) $\textit{, \ }$i,k=1,2,...n$\textit{\ in 
}$K_{n}$\textit{\ are defined by the relations } 
\begin{equation}
\sum\limits_{k=1}^{k=n}g^{ik}\left( \mathcal{P}^{n}\right) g_{lk}\left( 
\mathcal{P}^{n}\right) =\delta _{l}^{i},\qquad g_{il}\left( \mathcal{P}%
^{n}\right) =\left( \mathbf{P}_{0}\mathbf{P}_{i}.\mathbf{P}_{0}\mathbf{P}%
_{l}\right) ,\qquad i,l=1,2,...n  \label{a1.5b}
\end{equation}
\begin{equation}
F_{n}\left( \mathcal{P}^{n}\right) =\det \left\vert \left\vert g_{ik}\left( 
\mathcal{P}^{n}\right) \right\vert \right\vert \neq 0,\qquad i,k=1,2,...n
\label{g2.6}
\end{equation}

\textit{III. Linear structure of the Euclidean space: } 
\begin{equation}
\sigma \left( P,Q\right) =\frac{1}{2}\sum\limits_{i,k=1}^{i,k=n}g^{ik}\left( 
\mathcal{P}^{n}\right) \left( x_{i}\left( P\right) -x_{i}\left( Q\right)
\right) \left( x_{k}\left( P\right) -x_{k}\left( Q\right) \right) ,\qquad
\forall P,Q\in \Omega  \label{a1.5a}
\end{equation}
\textit{where coordinates }$x_{i}\left( P\right) ,$\textit{\ }$i=1,2,...n$%
\textit{\ of the point }$P$\textit{\ are covariant coordinates of the vector 
}$\mathbf{P}_{0}\mathbf{P}$\textit{, defined by the relation } 
\begin{equation}
x_{i}\left( P\right) =\left( \mathbf{P}_{0}\mathbf{P}_{i}.\mathbf{P}_{0}%
\mathbf{P}\right) ,\qquad i=1,2,...n  \label{b12}
\end{equation}

\textit{IV: The metric tensor matrix }$g_{lk}\left( \mathcal{P}^{n}\right) $%
\textit{\ has only positive eigenvalues } 
\begin{equation}
g_{k}>0,\qquad k=1,2,...,n  \label{a15c}
\end{equation}

\textit{V. The continuity condition: the system of equations } 
\begin{equation}
\left( \mathbf{P}_{0}\mathbf{P}_{i}.\mathbf{P}_{0}\mathbf{P}\right)
=y_{i}\in \mathbb{R},\qquad i=1,2,...n  \label{b14}
\end{equation}
\textit{considered to be equations for determination of the point }$P$%
\textit{\ as a function of coordinates }$y=\left\{ y_{i}\right\} $\textit{,\
\ }$i=1,2,...n$\textit{\ has always one and only one solution. }Conditions
II -- V contain a reference to the dimension $n$\ of the Euclidean space.

This theorem states that the proper Euclidean space has the property of the $%
\sigma $-immanence, and hence any statement $\mathcal{S}$ of the proper
Euclidean geometry can be expressed in terms and only in terms of the world
function $\sigma _{\mathrm{E}}$ of the Euclidean geometry in the form $%
\mathcal{S}\left( \sigma _{\mathrm{E}}\right) $. Producing the change $%
\sigma _{\mathrm{E}}\rightarrow \sigma $ in the statement $\mathcal{S}$, we
obtain corresponding statement $\mathcal{S}\left( \sigma \right) $ of
another T-geometry $\mathcal{G}$, described by the world function $\sigma $.

\section{Construction of geometric objects in the \newline
$\protect\sigma $-immanent form}

In the T-geometry the geometric object $\mathcal{O}$ is described by means
of the skeleton-envelope method \cite{R01}. It means that any geometric
object $\mathcal{O}$ is considered to be a set of intersections and joins of
elementary geometric objects (EGO).

The finite set $\mathcal{P}^{n}\equiv \left\{ P_{0},P_{1},...,P_{n}\right\}
\subset \Omega $ of parameters of the envelope function $f_{\mathcal{P}^{n}}$
is the skeleton of elementary geometric object (EGO) $\mathcal{E}\subset
\Omega $. The set $\mathcal{E}\subset \Omega $ of points forming EGO is
called the envelope of its skeleton $\mathcal{P}^{n}$. In the continuous
generalized geometry the envelope $\mathcal{E}$ is usually a continual set
of points. The envelope function $f_{\mathcal{P}^{n}}$%
\begin{equation}
f_{\mathcal{P}^{n}}:\qquad \Omega \rightarrow \mathbb{R},  \label{h2.1}
\end{equation}
determining EGO is a function of the running point $R\in \Omega $ and of
parameters $\mathcal{P}^{n}\subset \Omega $. The envelope function $f_{%
\mathcal{P}^{n}}$ is supposed to be an algebraic function of $s$ arguments $%
w=\left\{ w_{1},w_{2},...w_{s}\right\} $, $s=(n+2)(n+1)/2$. Each of
arguments $w_{k}=\sigma \left( Q_{k},L_{k}\right) $ is a $\sigma $-function
of two arguments $Q_{k},L_{k}\in \left\{ R,\mathcal{P}^{n}\right\} $, either
belonging to skeleton $\mathcal{P}^{n}$, or coinciding with the running
point $R$. Thus, any elementary geometric object $\mathcal{E}$ is determined
by its skeleton and its envelope function as the set of zeros of the
envelope function 
\begin{equation}
\mathcal{E}=\left\{ R|f_{\mathcal{P}^{n}}\left( R\right) =0\right\}
\label{h2.2}
\end{equation}

For instance, the cylinder $\mathcal{C}(P_{0},P_{1},Q)$ with the points $%
P_{0},P_{1}$ on the cylinder axis and the point $Q$ on its surface is
determined by the relation 
\begin{eqnarray}
\mathcal{C}(P_{0},P_{1},Q) &=&\left\{ R|f_{P_{0}P_{1}Q}\left( R\right)
=0\right\} ,  \label{g3.1} \\
f_{P_{0}P_{1}Q}\left( R\right) &=&F_{2}\left( P_{0},P_{1},Q\right)
-F_{2}\left( P_{0},P_{1},R\right)  \nonumber
\end{eqnarray}
\begin{equation}
F_{2}\left( P_{0},P_{1},Q\right) =\left\vert 
\begin{array}{cc}
\left( \mathbf{P}_{0}\mathbf{P}_{1}.\mathbf{P}_{0}\mathbf{P}_{1}\right) & 
\left( \mathbf{P}_{0}\mathbf{P}_{1}.\mathbf{P}_{0}\mathbf{Q}\right) \\ 
\left( \mathbf{P}_{0}\mathbf{Q}.\mathbf{P}_{0}\mathbf{P}_{1}\right) & \left( 
\mathbf{P}_{0}\mathbf{Q}.\mathbf{P}_{0}\mathbf{Q}\right)%
\end{array}
\right\vert  \label{g3.2}
\end{equation}
Here $\frac{1}{2}\sqrt{F_{2}\left( P_{0},P_{1},Q\right) }\ $ is the area of
triangle with vertices at the points $P_{0},P_{1},Q$. The equality $%
F_{2}\left( P_{0},P_{1},Q\right) =F_{2}\left( P_{0},P_{1},R\right) $ means
that the distance between the point $Q$ and the axis, determined by the
vector $\mathbf{P}_{0}\mathbf{P}_{1}$ is equal to the distance between $R$
and the axis.

The elementary geometrical object $\mathcal{E}$ is determined in all
T-geometries at once. In particular, it is determined in the proper
Euclidean geometry, where we can obtain its meaning. We interpret the
elementary geometrical object $\mathcal{E}$, using our knowledge of the
proper Euclidean geometry. Thus, the proper Euclidean geometry is used as a
sample geometry for interpretation of any generalized geometry. In
particular, the cylinder (\ref{g3.1}) is determined uniquely in any
T-geometry with any world function $\sigma .$

In the Euclidean geometry the points $P_{0}$ and $P_{1}$ determine the
cylinder axis. The shape of a cylinder depends on its axis and radius, but
not on the disposition of points $P_{0},P_{1}$ on the cylinder axis. As a
result in the Euclidean geometry the cylinders $\mathcal{C}(P_{0},P_{1},Q)$
and $\mathcal{C}(P_{0},P_{2},Q)$ coincide, provided vectors $\mathbf{P}_{0}%
\mathbf{P}_{1}$ and $\mathbf{P}_{0}\mathbf{P}_{2}$ are collinear. In the
general case of T-geometry the cylinders $\mathcal{C}(P_{0},P_{1},Q)$ and $%
\mathcal{C}(P_{0},P_{2},Q)$ do not coincide, in general, even if vectors $%
\mathbf{P}_{0}\mathbf{P}_{1}$ and $\mathbf{P}_{0}\mathbf{P}_{2}$ are
collinear. Thus, in general, a deformation of the Euclidean geometry splits
the Euclidean geometrical objects.

At construction of a generalized geometry we do not try to repeat derivation
of Euclidean algorithms from other axioms. We take the geometrical objects
and relations between them, prepared in the framework of the Euclidean
geometry and describe them in terms of the world function. Thereafter we
deform them, replacing the Euclidean world function $\sigma _{\mathrm{E}}$
by the world function $\sigma $ of the geometry in question. In practice the
construction of the elementary geometric object is reduced to the
representation of the corresponding Euclidean geometrical object in the $%
\sigma $-immanent form, i.e. in terms of the Euclidean world function. The
last problem is the problem of the proper Euclidean geometry. The problem of
representation of the geometrical object (or relation between objects) in
the $\sigma $-immanent form is a real problem of the T-geometry construction.

Application of the deformation principle is restricted by two constraints.

1. The deformation principle is to be applied separately from other methods
of the geometry construction. In particular, one may not use topological
structures in construction of a T-geometry, because for effective
application of the deformation principle the obtained T-geometry must be
determined only by the world function.

2. Describing Euclidean geometric objects $\mathcal{O}\left( \sigma _{%
\mathrm{E}}\right) $ and Euclidean relation $\mathcal{R}\left( \sigma _{%
\mathrm{E}}\right) $ in terms of $\sigma _{\mathrm{E}}$, we are not to use
special properties of Euclidean world function $\sigma _{\mathrm{E}}$. In
particular, definitions of $\mathcal{O}\left( \sigma _{\mathrm{E}}\right) $
and $\mathcal{R}\left( \sigma _{\mathrm{E}}\right) $ are to have similar
form in Euclidean geometries of different dimensions. They must not depend
on the dimension of the Euclidean space.

The T-geometry construction is not to use coordinates and other methods of
description, because the application of the means of description imposes
constraints on the constructed geometry. Any means of description is a
structure $St$ given on the basic Euclidean geometry with the world function 
$\sigma _{\mathrm{E}}$. Replacement $\sigma _{\mathrm{E}}\rightarrow \sigma $
is sufficient for construction of unique generalized geometry $\mathcal{G}%
_{\sigma }$. If we use an additional structure $St$ for the T-geometry
construction, we obtain, in general, other geometry $\mathcal{G}_{St}$,
which coincides with $\mathcal{G}_{\sigma }$ not for all $\sigma $, but only
for some of world functions $\sigma $. Thus, a use of additional means of
description restricts the list of possible generalized geometries. For
instance, if we use the coordinate description at construction of the
generalized geometry, the obtained geometry appears to be continuous,
because description by means of the coordinates is effective only for
continuous geometries, where the number of coordinates coincides with the
geometry dimension.

As far as the $\sigma $-immanent description of the proper Euclidean
geometry is possible, it is possible for any T-geometry, because any
geometrical object $\mathcal{O}$ and any relation $\mathcal{R}$ in the
generalized geometry $\mathcal{G}$ is obtained from the corresponding
geometrical object $\mathcal{O}_{\mathrm{E}}$ and from the corresponding
relation $\mathcal{R}_{\mathrm{E}}$ in the proper Euclidean geometry $%
\mathcal{G}_{\mathrm{E}}$ by means of the replacement $\sigma _{\mathrm{E}%
}\rightarrow \sigma $ in description of $\mathcal{O}_{\mathrm{E}}$ and $%
\mathcal{R}_{\mathrm{E}}$. For such a replacement be possible, the
description of $\mathcal{O}_{\mathrm{E}}$ and $\mathcal{R}_{\mathrm{E}}$ is
not to refer to special properties of $\sigma _{\mathrm{E}}$, described by
conditions II -- V. A formal indicator of the conditions II -- V application
is a reference to the dimension $n$, because any of conditions II -- V
contains a reference to the dimension $n$ of the proper Euclidean space.

Let us suppose that some geometrical object $\mathcal{O}_{\mathrm{E}%
_{n}}\left( \sigma _{\mathrm{E}_{n}},n\right) $ is defined in the $n$%
-dimensional Euclidean space, and this definition refers explicitly to the
dimension of the Euclidean space $n$. We deform the $n$-dimensional
Euclidean space $E_{n}$ in the $m$-dimensional Euclidean space $E_{m}$. Then
we must make the change 
\begin{equation}
\mathcal{O}_{\mathrm{E}_{n}}\left( \sigma _{\mathrm{E}_{n}},n\right)
\rightarrow \mathcal{O}_{\mathrm{E}_{n}}\left( \sigma _{\mathrm{E}%
_{m}},n\right)  \label{g3.3}
\end{equation}
On the other hand, we may define the same geometrical object directly in the 
$m$-dimensional Euclidean space $E_{m}$ in the form $\mathcal{O}_{\mathrm{E}%
_{m}}\left( \sigma _{\mathrm{E}_{m}},m\right) $. Equating this expression to
(\ref{g3.3}), we obtain 
\begin{equation}
\mathcal{O}_{\mathrm{E}_{n}}\left( \sigma _{\mathrm{E}_{m}},n\right) =%
\mathcal{O}_{\mathrm{E}_{m}}\left( \sigma _{\mathrm{E}_{m}},m\right) ,\qquad
\forall m,n\in \mathbb{N}  \label{g3.4}
\end{equation}

It means that the definition of the geometrical object $\mathcal{O}$ is to
be independent on the dimension of the Euclidean space.

If nevertheless we use one of special properties II -- V of the Euclidean
space in the $\sigma $-immanent description of a geometrical object $%
\mathcal{O}$, or relation $\mathcal{R}$ , we refer to the dimension $n$ and,
ultimately, to the coordinate system, which is only a means of description.

Let us show this in the example of the determination of the straight in the
Euclidean space. The straight $\mathcal{T}_{P_{0}Q}$ in the proper Euclidean
space is defined by two its points $P_{0}$ and $Q$ $\;\left( P_{0}\neq
Q\right) $ as the set of points $R$ 
\begin{equation}
\mathcal{T}_{P_{0}Q}=\left\{ R\;|\;\mathbf{P}_{0}\mathbf{Q}||\mathbf{P}_{0}%
\mathbf{R}\right\}  \label{b15}
\end{equation}
where condition $\mathbf{P}_{0}\mathbf{Q}||\mathbf{P}_{0}\mathbf{R}$ means
that vectors $\mathbf{P}_{0}\mathbf{Q}$ and $\mathbf{P}_{0}\mathbf{R}$ are
collinear, i.e. the scalar product $\left( \mathbf{P}_{0}\mathbf{Q}.\mathbf{P%
}_{0}\mathbf{R}\right) $ of these two vectors satisfies the relation 
\begin{equation}
\mathbf{P}_{0}\mathbf{Q}||\mathbf{P}_{0}\mathbf{R:\qquad }\left( \mathbf{P}%
_{0}\mathbf{Q}.\mathbf{P}_{0}\mathbf{R}\right) ^{2}=\left( \mathbf{P}_{0}%
\mathbf{Q}.\mathbf{P}_{0}\mathbf{Q}\right) \left( \mathbf{P}_{0}\mathbf{R}.%
\mathbf{P}_{0}\mathbf{R}\right)  \label{b16}
\end{equation}
where the scalar product is defined by the relation (\ref{b1.1a}). Thus, the
straight line $\mathcal{T}_{P_{0}Q}$ is defined $\sigma $-immanently, i.e.
in terms of the world function $\sigma $. We shall use two different names
(straight and tube) for the geometric object $\mathcal{T}_{P_{0}Q}$. We
shall use the term \textquotedblright straight\textquotedblright , when we
want to stress that $\mathcal{T}_{P_{0}Q}$ is a result of deformation of the
Euclidean straight. We shall use the term \textquotedblright
tube\textquotedblright , when we want to stress that $\mathcal{T}_{P_{0}Q}$
may be a many-dimensional surface.

In the Euclidean geometry one can use another definition of collinearity.
Vectors $\mathbf{P}_{0}\mathbf{Q}$ and $\mathbf{P}_{0}\mathbf{R}$ are
collinear, if components of vectors $\mathbf{P}_{0}\mathbf{Q}$ and $\mathbf{P%
}_{0}\mathbf{R}$ are proportional in some rectilinear coordinate system. For
instance, in the $n$-dimensional Euclidean space one can introduce
rectilinear coordinate system, choosing $n+1$ points $\mathcal{P}%
^{n}=\left\{ P_{0},P_{1},...P_{n}\right\} $ and forming $n$ basic vectors $%
\mathbf{P}_{0}\mathbf{P}_{i}$, $i=1,2,...n$. Then the collinearity condition
can be written in the form of $n$ equations 
\begin{equation}
\mathbf{P}_{0}\mathbf{Q}||\mathbf{P}_{0}\mathbf{R:\qquad }\left( \mathbf{P}%
_{0}\mathbf{P}_{i}.\mathbf{P}_{0}\mathbf{Q}\right) =a\left( \mathbf{P}_{0}%
\mathbf{P}_{i}.\mathbf{P}_{0}\mathbf{R}\right) ,\qquad i=1,2,...n,\qquad
a\in \mathbb{R}\backslash \left\{ 0\right\}  \label{b17}
\end{equation}
where $a\neq 0$ is some real constant. Relations (\ref{b17}) are relations
for covariant components of vectors $\mathbf{P}_{0}\mathbf{Q}$ and $\mathbf{P%
}_{0}\mathbf{R}$ in the considered coordinate system with basic vectors $%
\mathbf{P}_{0}\mathbf{P}_{i}$, $i=1,2,...n$. The definition of collinearity (%
\ref{b17}) depends on the dimension $n$ of the Euclidean space. Let points $%
\mathcal{P}^{n}$ be chosen in such a way, that $\left( \mathbf{P}_{0}\mathbf{%
P}_{1}.\mathbf{P}_{0}\mathbf{Q}\right) \neq 0$. Then eliminating the
parameter $a$ from relations (\ref{b17}), we obtain $n-1$ independent
relations, and the geometrical object 
\begin{eqnarray}
\mathcal{T}_{Q\mathcal{P}^{n}} &=&\left\{ R\;|\;\mathbf{P}_{0}\mathbf{Q}||%
\mathbf{P}_{0}\mathbf{R}\right\} =\bigcap\limits_{i=2}^{i=n}\mathcal{S}_{i},
\label{c2.1} \\
\mathcal{S}_{i} &=&\left\{ R\left\vert \frac{\left( \mathbf{P}_{0}\mathbf{P}%
_{i}.\mathbf{P}_{0}\mathbf{Q}\right) }{\left( \mathbf{P}_{0}\mathbf{P}_{1}.%
\mathbf{P}_{0}\mathbf{Q}\right) }=\frac{\left( \mathbf{P}_{0}\mathbf{P}_{i}.%
\mathbf{P}_{0}\mathbf{R}\right) }{\left( \mathbf{P}_{0}\mathbf{P}_{1}.%
\mathbf{P}_{0}\mathbf{R}\right) }\right. \right\} ,\qquad i=2,3,...n
\label{c2.2}
\end{eqnarray}
defined according to (\ref{b17}), depends on $n+2$ points $Q,\mathcal{P}^{n}$%
. This geometrical object $\mathcal{T}_{Q\mathcal{P}^{n}}$ is defined $%
\sigma $-immanently. It is a complex, consisting of the straight line and of
the coordinate system, represented by $n+1$ points $\mathcal{P}^{n}=\left\{
P_{0},P_{1},...P_{n}\right\} $. In the Euclidean space the dependence on the
choice of the coordinate system and on $n$ points $\left\{
P_{1},...P_{n}\right\} $, determining this system, is fictitious. The
geometrical object $\mathcal{T}_{Q\mathcal{P}^{n}}$ depends essentially only
on two points $P_{0},Q$ and coincides with the straight line $\mathcal{T}%
_{P_{0}Q}$ in the Euclidean space. But at deformations of the Euclidean
space the geometrical objects $\mathcal{T}_{Q\mathcal{P}^{n}}$ and $\mathcal{%
T}_{P_{0}Q}$ are deformed differently. The points $P_{1},P_{2},...P_{n}$
cease to be fictitious in definition of $\mathcal{T}_{Q\mathcal{P}^{n}}$,
and geometrical objects $\mathcal{T}_{Q\mathcal{P}^{n}}$ and $\mathcal{T}%
_{P_{0}Q}$ become to be different geometric objects, in general. But being
different, in general, they may coincide in some special cases.

What of the two geometrical objects in the deformed geometry $\mathcal{G}$
should be interpreted as a straight line, passing through the points $P_{0}$
and $Q$ in the geometry $\mathcal{G}$? Of course, it is $\mathcal{T}%
_{P_{0}Q} $, because its definition does not contain a reference to a
coordinate system, whereas definition of $\mathcal{T}_{Q\mathcal{P}^{n}}$
depends on the choice of the coordinate system, represented by points $%
\mathcal{P}^{n}$. In general, definitions of geometric objects and relations
between them are not to refer to the means of description. Otherwise, the
points determining the coordinate system are to be included in definition of
the geometrical object.

But in the given case the geometrical object $\mathcal{T}_{P_{0}Q}$ is a $%
(n-1)$-dimensional surface, in general, whereas $\mathcal{T}_{Q\mathcal{P}%
^{n}}$ is an intersection of $(n-1)\;\;$ $(n-1)$-dimensional surfaces, i.e. $%
\mathcal{T}_{Q\mathcal{P}^{n}}$ is a one-dimensional curve, in general. The
one-dimensional curve $\mathcal{T}_{Q\mathcal{P}^{n}}$ corresponds better to
our ideas on the straight line, than the $(n-1)$-dimensional surface $%
\mathcal{T}_{P_{0}Q}$. Nevertheless, in T-geometry $\mathcal{G}$ it is $%
\mathcal{T}_{P_{0}Q}$, that is an analog of the Euclidean straight line.

It is very difficult to overcome our conventional idea that the Euclidean
straight line cannot be deformed into many-dimensional surface, and \textit{%
this idea has been prevent for years from construction of T-geometries}.
Practically one uses such generalized geometries, where deformation of the
Euclidean space transforms the Euclidean straight lines into one-dimensional
lines. It means that one chooses such geometries, where geometrical objects $%
\mathcal{T}_{P_{0}Q}$ and $\mathcal{T}_{Q\mathcal{P}^{n}}$ coincide. 
\begin{equation}
\mathcal{T}_{P_{0}Q}=\mathcal{T}_{Q\mathcal{P}^{n}}  \label{b19}
\end{equation}
Condition (\ref{b19}) of coincidence of the objects $\mathcal{T}_{P_{0}Q}$
and $\mathcal{T}_{Q\mathcal{P}^{n}}$, imposed on the T-geometry, restricts
the list of possible T-geometries.

In general, the condition (\ref{b19}) cannot be fulfilled, because lhs does
not depend on points $\left\{ P_{1},P_{2},...P_{n}\right\} $, whereas rhs of
(\ref{b19}) depends, in general. The tube $\mathcal{T}_{Q\mathcal{P}^{n}}$
does not depend on the points $\left\{ P_{1},P_{2},...P_{n}\right\} $,
provided the distance $\sqrt{2\sigma \left( P_{i},P_{k}\right) }$ between
any two points $P_{i},P_{k}\in \mathcal{P}^{n}$ is infinitesimal. In the
Riemannian geometry the constraint (\ref{b19}) is fulfilled at the
additional restriction. 
\begin{equation}
\sqrt{2\sigma \left( P_{i},P_{k}\right) }=\text{infinitesimal},\qquad
i,k=1,2,...n  \label{b20}
\end{equation}

\section{Interplay between metric geometry and \newline
T-geometry}

Let us consider the metric geometry, given on the set $\Omega $ of points.
The metric space $M=\left\{ \rho ,\Omega \right\} $ is given by the metric
(distance) $\rho $. 
\begin{eqnarray}
\rho &:&\quad \Omega \times \Omega \rightarrow \lbrack 0,\infty )\subset 
\mathbb{R}  \label{c2.3} \\
\rho (P,P) &=&0,\qquad \rho (P,Q)=\rho (Q,P),\qquad \forall P,Q\in \Omega
\label{c2.4} \\
\rho (P,Q) &\geq &0,\qquad \rho (P,Q)=0,\quad \text{iff }P=Q,\qquad \forall
P,Q\in \Omega  \label{c2.5} \\
0 &\leq &\rho (P,R)+\rho (R,Q)-\rho (P,Q),\qquad \forall P,Q,R\in \Omega
\label{c2.6}
\end{eqnarray}
At first sight the metric space is a special case of the $\sigma $-space (%
\ref{h1.1}), and the metric geometry is a special case of the T-geometry
with additional constraints (\ref{c2.5}), (\ref{c2.6}) imposed on the world
function $\sigma =\frac{1}{2}\rho ^{2}$. However it is not so, because the
metric geometry is not equipped by the deformation principle. The metric
geometry does not use the algorithms of the Euclidean geometry construction.
In the metric geometry the deformation principle can be used only in its
coordinate form 
\begin{equation}
g_{\mathrm{E}ik}dx^{i}dx^{k}\rightarrow g_{ik}dx^{i}dx^{k}  \label{c2.8a}
\end{equation}
because the coordinateless form (\ref{h1.2}) of the deformation principle as
well as the $\sigma $-immanence of the Euclidean geometry and complex of
conditions (\ref{a1.4}) - (\ref{b14}) were not known until 1990, although
each of relations (\ref{a1.4}) - (\ref{b14}) was well known. But the metric
geometry is described in the coordinateless form, and application of the
deformation (\ref{c2.8a}) is impossible in the metric geometry.

Additional (with respect to the $\sigma $-space) constraints (\ref{c2.5}), (%
\ref{c2.6}) are imposed to provide one-dimensionality of the straight lines.
In the metric geometry the shortest (straight) line can be constructed only
in the case, when it is one-dimensional.

Let us consider the set $\mathcal{EL}\left( P,Q,a\right) $ of points $R$ 
\begin{equation}
\mathcal{EL}\left( P,Q,2a\right) =\left\{ R|f_{P,Q,2a}\left( R\right)
=0\right\} ,\qquad f_{P,Q,2a}\left( R\right) =\rho (P,R)+\rho (R,Q)-2a
\label{c2.8}
\end{equation}
If the metric space coincides with the proper Euclidean space, this set of
points is an ellipsoid with focuses at the points $P,Q$ and the large
semiaxis $a$. The relations $f_{P,Q,2a}\left( R\right) >0$, $%
f_{P,Q,2a}\left( R\right) =0$, $f_{P,Q,2a}\left( R\right) <0$ determine
respectively external points, boundary points and internal points of the
ellipsoid. If $\rho \left( P,Q\right) =2a$, we obtain the degenerate
ellipsoid, which coincides with the segment $\mathcal{T}_{\left[ PQ\right] }$
of the straight line, passing through the points $P$, $Q$.\ In the proper
Euclidean geometry, the degenerate ellipsoid is one-dimensional segment of
the straight line, but it is not evident that it is one-dimensional in the
case of arbitrary metric geometry. For such a degenerate ellipsoid be
one-dimensional in the arbitrary metric space, it is necessary that any
degenerate ellipsoid $\mathcal{EL}\left( P,Q,\rho \left( P,Q\right) \right) $
have no internal points. This constraint is written in the form 
\begin{equation}
f_{P,Q,\rho \left( P,Q\right) }\left( R\right) =\rho (P,R)+\rho (R,Q)-\rho
(P,Q)\geq 0  \label{c2.9}
\end{equation}

Comparing relation (\ref{c2.9}) with (\ref{c2.6}), we see that the
constraint (\ref{c2.6}) is introduced to provide the straight (shortest)
line one-dimensionality (absence of internal points in the geometrical
object determined by two points).

As far as the metric geometry does not use the deformation principle, it is
a poor geometry, because in the framework of this geometry one cannot
construct the scalar product of two vectors, define linear independence of
vectors and construct such geometrical objects as planes. All these objects
as well as others are constructed on the basis of the deformation of the
proper Euclidean geometry.

Generalizing the metric geometry, Menger \cite{M28} and Blumenthal \cite{B53}
removed the triangle axiom (\ref{c2.6}). They tried to construct the
distance geometry, which would be a more general geometry, than the metric
one. As far as they did not use the deformation principle, they could not
determine the shortest (straight) line without a reference to the
topological concept of the curve $\mathcal{L}$, defined as a continuous
mapping (\ref{g1.1}), which cannot be expressed only via the distance. As a
result the distance geometry appeared to be not a pure metric geometry (i.e.
the geometry determined only by the distance).

Note that the Riemannian geometry uses the deformation principle in the
coordinate form. The distance geometry cannot use it in such a form, because
the metric and distance geometries are formulated in the coordinateless
form. It is to use the deformation principle in the coordinateless form. But
application of the deformation principle in the coordinateless form needs a
use of the Euclidean geometry $\sigma $-immanence. K. Menger went to the
concept of the $\sigma $-immanence, but he stopped in one step before the $%
\sigma $-immanence. Look at the Menger's theorem \cite{M28}, written in our
designations

\begin{theorem}
The $\sigma $-space $V=\left\{ \sigma ,\Omega \right\} $ is
isometrically embeddable in $n$-dimensional proper Euclidean space
$E_{n}$, if and only if any set of $n+3$ points of $\Omega $ is
isometrically embeddable in $E_{n}$.
\end{theorem}

The theorem on the $\sigma $-immanence of the Euclidean geometry is obtained
from the Menger's theorem, if instead of the condition "any set of $n+3$
points of $\Omega $ is isometrically embeddable in $E_{n}"$ one writes the
condition (\ref{a1.5a}), which also contains $n+3$ points: $P,Q,\mathcal{P}%
^{n}$ and describes the fact that $\left\{ P,Q,\mathcal{P}^{n}\right\}
\subset E_{n}$. In this case the theorem condition contains only a reference
to the properties of the world function of the Euclidean space, but not to
the Euclideaness of the space. (continuity of the $\sigma $-space $V$ is
neglected in such a formulation.)

\section{Conditions of the deformation principle \newline
application}

Riemannian geometries satisfy the condition (\ref{b19}). The Riemannian
geometry is a kind of inhomogeneous generalized geometry, and, hence, it
uses the deformation principle. Constructing the Riemannian geometry, the
infinitesimal Euclidean distance is deformed into the Riemannian distance.
The deformation is chosen in such a way that any Euclidean straight line $%
\mathcal{T}_{\mathrm{E}P_{0}Q}$, passing through the point $P_{0}$,
collinear to the vector $\mathbf{P}_{0}\mathbf{Q}$, is transformed into the
geodesic $\mathcal{T}_{P_{0}Q}$, passing through the point $P_{0}$,
collinear to the vector $\mathbf{P}_{0}\mathbf{Q}$ in the Riemannian space.

Note that in T-geometries, satisfying the condition (\ref{b19}) for all
points $Q,\mathcal{P}^{n}$, the straight line 
\begin{equation}
\mathcal{T}_{Q_{0};P_{0}Q}=\left\{ R\;|\;\mathbf{P}_{0}\mathbf{Q}||\mathbf{Q}%
_{0}\mathbf{R}\right\}  \label{b3.0}
\end{equation}
passing through the point $Q_{0}$ collinear to the vector $\mathbf{P}_{0}%
\mathbf{Q}$, is not a one-dimensional line, in general. If the Riemannian
geometries be T-geometries, they would contain non-one-dimensional geodesics
(straight lines). But the Riemannian geometries are not T-geometries,
because at their construction one uses not only the deformation principle,
but some other methods, containing a reference to the means of description.
In particular, in the Riemannian geometries the absolute parallelism is
absent, and one cannot define a straight line (\ref{b3.0}), because the
collinearity relation $\mathbf{P}_{0}\mathbf{Q}||\mathbf{Q}_{0}\mathbf{R}$
is not defined, if points $P_{0}$ and $Q_{0}$ do not coincide. On one hand,
a lack of absolute parallelism allows one to go around the problem of
non-one-dimensional straight lines. On the other hand, it makes the
Riemannian geometries to be inconsistent, because they cease to be
T-geometries, which are consistent by the construction (see for details \cite%
{R02}).

The fact is that the application of \textit{only deformation principle }is
sufficient for construction of a generalized geometry. Besides, such a
construction is consistent, because the original Euclidean geometry is
consistent and, deforming it, we do not use any logical reasonings. If we
introduce additional structure (for instance, a topological structure) we
obtain a fortified geometry, i.e. a generalized geometry with additional
structure on it. The T-geometry, equipped with additional structure, is a
more pithy construction, than the T-geometry simply. But it is valid only in
the case, when we consider the additional structure as an addition to the
T-geometry. If we use an additional structure in construction of the
geometry, we identify the additional structure with one of structures of the
geometry. If we demand that the additional structure be a structure of
T-geometry, we restrict an application of the deformation principle and
reduce the list of possible generalized geometries, because coincidence of
the additional structure with some structure of a geometry is possible not
for all geometries, but only for some of them.

Let, for instance, we use concept of a curve $\mathcal{L}$ (\ref{g1.1}) for
construction of a generalized geometry. The concept of curve $\mathcal{L}$,
considered as a continuous mapping, is a topological structure, which cannot
be expressed only via the distance or via the world function. A use of the
mapping (\ref{g1.1}) needs an introduction of topological space and, in
particular, the concept of continuity. If we identify the topological curve (%
\ref{g1.1}) with the ''metrical''\ curve, defined as a broken line 
\begin{equation}
\mathcal{T}_{\mathrm{br}}=\bigcup\limits_{i}\mathcal{T}_{\left[ P_{i}P_{i+1}%
\right] },\qquad \mathcal{T}_{\left[ P_{i}P_{i+1}\right] }=\left\{ R|\sqrt{%
2\sigma \left( P_{i},P_{i+1}\right) }-\sqrt{2\sigma \left( P_{i},R\right) }-%
\sqrt{2\sigma \left( R,P_{i+1}\right) }\right\}  \label{a1.2}
\end{equation}
consisting of the straight line segments $\mathcal{T}_{\left[ P_{i}P_{i+1}%
\right] }$ between the points $P_{i}$, $P_{i+1}$, we truncate the list of
possible geometries, because such an identification is possible only in some
generalized geometries. Identifying (\ref{g1.1}) and (\ref{a1.2}), we
eliminate all discrete geometries and those continuous geometries, where the
segment $\mathcal{T}_{\left[ P_{i}P_{i+1}\right] }$ of straight line is a
surface, but not a one-dimensional set of points. Thus, additional
structures may lead to (i) a fortified geometry, (ii) a restricted geometry
and (iii) a restricted fortified geometry. The result depends on the method
of the additional structure application.

Note that some constraints (continuity, convexity, lack of absolute
parallelism), imposed on generalized geometries, are a result of a
disagreement of the means of description, which are used at the geometry
construction. In the T-geometry, which uses only the deformation principle,
there is no such restrictions. Besides, the T-geometry has some new property
of a geometry, which is not accepted by conventional versions of generalized
geometry. This property, called the geometry nondegeneracy, follows directly
from the application of arbitrary deformations to the proper Euclidean
geometry.

\begin{definition}
The geometry is degenerate at the point $P_{0}$ in the direction of the
vector $\mathbf{Q}_{0}\mathbf{Q}$, $\left\vert \mathbf{Q}_{0}\mathbf{Q}%
\right\vert \neq 0$, if the relations
\begin{equation}
\mathbf{Q}_{0}\mathbf{Q}\uparrow \uparrow \mathbf{P}_{0}\mathbf{R:\qquad }%
\left( \mathbf{Q}_{0}\mathbf{Q}.\mathbf{P}_{0}\mathbf{R}\right) =\sqrt{%
\left\vert \mathbf{Q}_{0}\mathbf{Q}\right\vert \cdot \left\vert \mathbf{P}%
_{0}\mathbf{R}\right\vert },\qquad \left\vert \mathbf{P}_{0}\mathbf{R}%
\right\vert =a\neq 0  \label{b3.1}
\end{equation}%
considered as equations for determination of the point $R$, have not more,
than one solution for any $a\neq 0$. Otherwise, the geometry is
nondegenerate at the point $P_{0}$ in the direction of the vector $\mathbf{Q}%
_{0}\mathbf{Q}$.
\end{definition}

\noindent Note that the first equation (\ref{b3.1}) is the condition of the
parallelism of vectors $\mathbf{Q}_{0}\mathbf{Q}$ and $\mathbf{P}_{0}\mathbf{%
R}$.

The proper Euclidean geometry is degenerate, i.e. it is degenerate at all
points in directions of all vectors. Considering the Minkowski geometry, one
should distinguish between the Minkowski T-geometry and Minkowski geometry.
The two geometries are described by the same world function and differ in
the definition of the parallelism. In the Minkowski T-geometry the
parallelism of two vectors $\mathbf{\mathbf{Q}_{0}\mathbf{Q}}$ and $\mathbf{%
\mathbf{P}_{0}\mathbf{R}}$ is defined by the first equation (\ref{b3.1}).
This definition is based on the deformation principle. In the $n$%
-dimensional Minkowski geometry ($n$-dimensional pseudo-Euclidean geometry
of index $1$) the parallelism is defined by the relation of the type of (\ref%
{b17}) 
\begin{equation}
\mathbf{Q}_{0}\mathbf{Q}\uparrow \uparrow \mathbf{P}_{0}\mathbf{R:\qquad }%
\left( \mathbf{P}_{0}\mathbf{P}_{i}.\mathbf{Q}_{0}\mathbf{Q}\right) =a\left( 
\mathbf{P}_{0}\mathbf{P}_{i}.\mathbf{P}_{0}\mathbf{R}\right) ,\qquad
i=1,2,...n,\qquad a>0  \label{b3.1a}
\end{equation}
where points $\mathcal{P}^{n}=\left\{ P_{0},P_{1},...P_{n}\right\} $
determine a rectilinear coordinate system with basic vectors $\mathbf{P}_{0}%
\mathbf{P}_{i}$, $i=1,2,..n$ in the $n$-dimensional Minkowski space.
Dependence of the definition (\ref{b3.1a}) on the points $\left(
P_{0},P_{1},...P_{n}\right) $ is fictitious, but dependence on the number $%
n+1$ of points $\mathcal{P}^{n}$ is essential. Thus, definition (\ref{b3.1a}%
) depends on the method of the geometry description.

The Minkowski T-geometry is degenerate at all points in direction of all
timelike vectors, and it is nondegenerate at all points in direction of all
spacelike vectors. The Minkowski geometry is degenerate at all points in
direction of all vectors. Conventionally one uses the Minkowski geometry,
ignoring the nondegeneracy in spacelike directions.

Considering the proper Riemannian geometry, one should distinguish between
the Riemannian T-geometry and the Riemannian geometry. The two geometries
are described by the same world function. They differ in the definition of
the parallelism. In the Riemannian T-geometry the parallelism of two vectors 
$\mathbf{\mathbf{Q}_{0}\mathbf{Q}}$ and $\mathbf{\mathbf{P}_{0}\mathbf{R}}$
is defined by (\ref{b3.1}). In the Riemannian geometry the parallelism of
two vectors $\mathbf{\mathbf{Q}_{0}\mathbf{Q}}$ and $\mathbf{\mathbf{P}_{0}%
\mathbf{R}}$ is defined only in the case, when the points $P_{0}$ and $Q_{0}$
coincide. Parallelism of remote vectors $\mathbf{\mathbf{Q}_{0}\mathbf{Q}}$
and $\mathbf{\mathbf{P}_{0}\mathbf{R}}$ is not defined, in general. This
fact is known as absence of absolute parallelism.

The proper Riemannian T-geometry is locally degenerate, i.e. it is
degenerate at all points $P_{0}$ in direction of all vectors $\mathbf{P}_{0}%
\mathbf{Q}$ with the origin at the point $P_{0}$. In the general case, when $%
P_{0}\neq Q_{0}$, the proper Riemannian T-geometry is nondegenerate, in
general. But it is degenerate locally as well as the proper Riemannian
geometry. The proper Riemannian geometry is degenerate, because it is
degenerate locally, whereas the illocal degeneracy is not defined in the
Riemannian geometry, because of the lack of absolute parallelism.
Conventionally one uses the Riemannian geometry (not Riemannian T-geometry)
and ignores the property of the nondegeneracy completely.

From the viewpoint of the conventional approach to the generalized geometry
the nondegeneracy is an undesirable property of a generalized geometry,
although from the logical viewpoint and from viewpoint of the deformation
principle the nondegeneracy is \textit{an inherent property of a generalized
geometry}. The illocal nondegeneracy is ejected from the proper Riemannian
geometry by denial of existence of the remote vector parallelism.
Nondegeneracy in the spacelike directions is ejected from the Minkowski
geometry by means of the redefinition of the two vectors parallelism. But
the nondegeneracy is an important property of the real space-time geometry.
To appreciate this, let us consider an example.

\section{Simple example of nondegenerate space-time geometry}

Let the space-time geometry $\mathcal{G}_{\mathrm{d}}$ be described by the
T-geometry, given on 4-dimensional manifold $\mathcal{M}_{1+3}$. The world
function $\sigma _{\mathrm{d}}$ is described by the relation 
\begin{equation}
\sigma _{\mathrm{d}}=\sigma _{\mathrm{M}}+D\left( \sigma _{\mathrm{M}%
}\right) =\left\{ 
\begin{array}{ll}
\sigma _{\mathrm{M}}+d & \text{if\ }\sigma _{0}<\sigma _{\mathrm{M}} \\ 
\left( 1+\frac{d}{\sigma _{0}}\right) \sigma _{\mathrm{M}} & \text{if\ }%
0\leq \sigma _{\mathrm{M}}\leq \sigma _{0} \\ 
\sigma _{\mathrm{M}} & \text{if\ }\sigma _{\mathrm{M}}<0%
\end{array}
\right.  \label{b3.3}
\end{equation}
where $d\geq 0$ and $\sigma _{0}>0$ are some constants. The quantity $\sigma
_{\mathrm{M}}$ is the world function in the Minkowski space-time geometry $%
\mathcal{G}_{\mathrm{M}}$. In the orthogonal rectilinear (inertial)
coordinate system $x=\left\{ t,\mathbf{x}\right\} $ the world function $%
\sigma _{\mathrm{M}}$ has the form 
\begin{equation}
\sigma _{\mathrm{M}}\left( x,x^{\prime }\right) =\frac{1}{2}\left(
c^{2}\left( t-t^{\prime }\right) ^{2}-\left( \mathbf{x}-\mathbf{x}^{\prime
}\right) ^{2}\right)  \label{b3.4}
\end{equation}
where $c$ is the speed of the light.

Let us compare the broken line (\ref{a1.2}) in Minkowski space-time geometry 
$\mathcal{G}_{\mathrm{M}}$ and in the distorted geometry $\mathcal{G}_{%
\mathrm{d}}$. We suppose that $\mathcal{T}_{\mathrm{br}}$ is timelike broken
line, and all links $\mathcal{T}_{\left[ P_{i}P_{i+1}\right] }$ of $\mathcal{%
T}_{\mathrm{br}}$ are timelike and have the same length 
\begin{equation}
\left\vert \mathbf{P}_{i}\mathbf{P}_{i+1}\right\vert _{\mathrm{d}}=\sqrt{%
2\sigma _{\mathrm{d}}\left( P_{i},P_{i+1}\right) }=\mu _{\mathrm{d}%
}>0,\qquad i=0,\pm 1,\pm 2,...  \label{b3.5}
\end{equation}
\begin{equation}
\left\vert \mathbf{P}_{i}\mathbf{P}_{i+1}\right\vert _{\mathrm{M}}=\sqrt{%
2\sigma _{\mathrm{M}}\left( P_{i},P_{i+1}\right) }=\mu _{\mathrm{M}%
}>0,\qquad i=0,\pm 1,\pm 2,...  \label{b3.5a}
\end{equation}
where indices \textquotedblright d\textquotedblright\ and \textquotedblright
M\textquotedblright\ mean that the quantity is calculated by means of $%
\sigma _{\mathrm{d}}$ and $\sigma _{\mathrm{M}}$ respectively. Vector $%
\mathbf{P}_{i}\mathbf{P}_{i+1}$ is regarded as the momentum of the particle
at the segment $\mathcal{T}_{\left[ P_{i}P_{i+1}\right] }$, and the quantity 
$\left\vert \mathbf{P}_{i}\mathbf{P}_{i+1}\right\vert =\mu $ is interpreted
as its (geometric) mass. It follows from definition (\ref{b1.1a}) and
relation (\ref{b3.3}), that for timelike vectors $\mathbf{P}_{i}\mathbf{P}%
_{i+1}$ with $\mu >\sqrt{2\sigma _{0}}$%
\begin{equation}
\left\vert \mathbf{P}_{i}\mathbf{P}_{i+1}\right\vert _{\mathrm{d}}^{2}=\mu _{%
\mathrm{d}}^{2}=\mu _{\mathrm{M}}^{2}+2d,\qquad \mu _{\mathrm{M}%
}^{2}>2\sigma _{0}  \label{b3.6}
\end{equation}
\begin{equation}
\left( \mathbf{P}_{i-1}\mathbf{P}_{i}.\mathbf{P}_{i}\mathbf{P}_{i+1}\right)
_{\mathrm{d}}=\left( \mathbf{P}_{i-1}\mathbf{P}_{i}.\mathbf{P}_{i}\mathbf{P}%
_{i+1}\right) _{\mathrm{M}}+d  \label{b3.7}
\end{equation}
Calculation of the shape of the segment $\mathcal{T}_{\left[ P_{0}P_{1}%
\right] }\left( \sigma _{\mathrm{d}}\right) $ in $\mathcal{G}_{\mathrm{d}}$
gives the relation 
\begin{equation}
r^{2}(\tau )=\left\{ 
\begin{array}{ll}
\tau ^{2}\mu _{\mathrm{d}}^{2}\frac{\left( 1-\frac{\tau d}{2\left( \sigma
_{0}+d\right) }\right) ^{2}}{\left( 1-\frac{2d}{\mu _{\mathrm{d}}^{2}}%
\right) }-\frac{\tau ^{2}\mu _{\mathrm{d}}^{2}\sigma _{0}}{\left( \sigma
_{0}+d\right) }, & 0<\tau <\frac{\sqrt{2(\sigma _{0}+d)}}{\mu _{\mathrm{d}}}
\\ 
\frac{3d}{2}+2d\left( \tau -1/2\right) ^{2}\left( 1-\frac{2d}{\mu _{\mathrm{d%
}}^{2}}\right) ^{-1}, & \frac{\sqrt{2(\sigma _{0}+d)}}{\mu _{\mathrm{d}}}%
<\tau <1-\frac{\sqrt{2(\sigma _{0}+d)}}{\mu _{\mathrm{d}}} \\ 
\left( 1-\tau \right) ^{2}\mu _{\mathrm{d}}^{2}\left[ \frac{\left( 1-\frac{%
\left( 1-\tau \right) d}{2\left( \sigma _{0}+d\right) }\right) ^{2}}{\left(
1-\frac{2d}{\mu _{\mathrm{d}}^{2}}\right) }-\frac{\sigma _{0}}{\left( \sigma
_{0}+d\right) }\right] , & 1-\frac{\sqrt{2(\sigma _{0}+d)}}{\mu _{\mathrm{d}}%
}<\tau <1%
\end{array}
\right. ,  \label{b3.7a}
\end{equation}
where $r\left( \tau \right) $ is the spatial radius of the segment $\mathcal{%
T}_{\left[ P_{0}P_{1}\right] }\left( \sigma _{\mathrm{d}}\right) $ in the
coordinate system, where points $P_{0}$ and $P_{1}$ have coordinates $%
P_{0}=\left\{ 0,0,0,0\right\} $, $P_{1}=\left\{ \mu _{\mathrm{d}%
},0,0,0\right\} $ and $\tau $ is a parameter along the segment $\mathcal{T}_{%
\left[ P_{0}P_{1}\right] }\left( \sigma _{\mathrm{d}}\right) $, ($\tau
\left( P_{0}\right) =0$, $\tau \left( P_{1}\right) =1$). One can see from (%
\ref{b3.7a}) that the characteristic value of the segment radius is equal to 
$\sqrt{d}$.

Let the broken tube $\mathcal{T}_{\mathrm{br}}$ describe the
\textquotedblright world tube\textquotedblright\ of a free particle. It
means by definition that any link $\mathbf{P}_{i-1}\mathbf{P}_{i}$ is
parallel to the adjacent link $\mathbf{P}_{i}\mathbf{P}_{i+1}$%
\begin{equation}
\mathbf{P}_{i-1}\mathbf{P}_{i}\uparrow \uparrow \mathbf{P}_{i}\mathbf{P}%
_{i+1}:\qquad \left( \mathbf{P}_{i-1}\mathbf{P}_{i}.\mathbf{P}_{i}\mathbf{P}%
_{i+1}\right) -\left\vert \mathbf{P}_{i-1}\mathbf{P}_{i}\right\vert \cdot
\left\vert \mathbf{P}_{i}\mathbf{P}_{i+1}\right\vert =0  \label{b3.8}
\end{equation}
Definition of parallelism is different in geometries $\mathcal{G}_{\mathrm{M}%
}$ and $\mathcal{G}_{\mathrm{d}}$. As a result links, which are parallel in
the geometry $\mathcal{G}_{\mathrm{M}}$, are not parallel in $\mathcal{G}_{%
\mathrm{d}}$ and vice versa.

Let $\mathcal{T}_{\mathrm{br}}\left( \sigma _{\mathrm{M}}\right) $ describe
the world line of a free particle in the geometry $\mathcal{G}_{\mathrm{M}}$%
. The angle $\vartheta _{\mathrm{M}}$ between the adjacent links in $%
\mathcal{G}_{\mathrm{M}}$ is defined by the relation 
\begin{equation}
\cosh \vartheta _{\mathrm{M}}=\frac{\left( \mathbf{P}_{-1}\mathbf{P}_{0}.%
\mathbf{P}_{0}\mathbf{P}_{1}\right) _{\mathrm{M}}}{\left| \mathbf{P}_{0}%
\mathbf{P}_{1}\right| _{\mathrm{M}}\cdot \left| \mathbf{P}_{-1}\mathbf{P}%
_{0}\right| _{\mathrm{M}}}=1  \label{b3.9}
\end{equation}
The angle $\vartheta _{\mathrm{M}}=0$, and the geometrical object $\mathcal{T%
}_{\mathrm{br}}\left( \sigma _{\mathrm{M}}\right) $ is a timelike straight
line on the manifold $\mathcal{M}_{1+3}$.

Let now $\mathcal{T}_{\mathrm{br}}\left( \sigma _{\mathrm{d}}\right) $
describe the world tube of a free particle in the geometry $\mathcal{G}_{%
\mathrm{d}}$. The angle $\vartheta _{\mathrm{d}}$ between the adjacent links
in $\mathcal{G}_{\mathrm{d}}$ is defined by the relation 
\begin{equation}
\cosh \vartheta _{\mathrm{d}}=\frac{\left( \mathbf{P}_{i-1}\mathbf{P}_{i}.%
\mathbf{P}_{i}\mathbf{P}_{i+1}\right) _{\mathrm{d}}}{\left\vert \mathbf{P}%
_{i}\mathbf{P}_{i+1}\right\vert _{\mathrm{d}}\cdot \left\vert \mathbf{P}%
_{i-1}\mathbf{P}_{i}\right\vert _{\mathrm{d}}}=1  \label{b3.10}
\end{equation}%
The angle $\vartheta _{\mathrm{d}}=0$ also. If we draw the broken tube $%
\mathcal{T}_{\mathrm{br}}\left( \sigma _{\mathrm{d}}\right) $ on the
manifold $\mathcal{M}_{1+3}$, using coordinates of basic points $P_{i}$ and
measure the angle $\vartheta _{\mathrm{dM}}$ between the adjacent links in
the Minkowski geometry $\mathcal{G}_{\mathrm{M}}$, we obtain from (\ref{b3.6}%
), (\ref{b3.7}) the following relation for the angle $\vartheta _{\mathrm{dM}%
}$  
\begin{equation}
\cosh \vartheta _{\mathrm{dM}}=\frac{\left( \mathbf{P}_{i-1}\mathbf{P}_{i}.%
\mathbf{P}_{i}\mathbf{P}_{i+1}\right) _{\mathrm{M}}}{\left\vert \mathbf{P}%
_{i}\mathbf{P}_{i+1}\right\vert _{\mathrm{M}}\cdot \left\vert \mathbf{P}%
_{i-1}\mathbf{P}_{i}\right\vert _{\mathrm{M}}}=\frac{\left( \mathbf{P}_{i-1}%
\mathbf{P}_{i}.\mathbf{P}_{i}\mathbf{P}_{i+1}\right) _{\mathrm{d}}-d}{%
\left\vert \mathbf{P}_{i}\mathbf{P}_{i+1}\right\vert _{\mathrm{d}}^{2}-2d}
\label{b3.11}
\end{equation}%
Substituting the value of $\left( \mathbf{P}_{i-1}\mathbf{P}_{i}.\mathbf{P}%
_{i}\mathbf{P}_{i+1}\right) _{\mathrm{d}}$, taken from (\ref{b3.10}), (\ref%
{b3.6}), we obtain 
\begin{equation}
\cosh \vartheta _{\mathrm{dM}}=\frac{\mu _{\mathrm{d}}^{d}-d}{\mu _{\mathrm{d%
}}^{2}-2d}\approx 1+\frac{d}{\mu _{\mathrm{d}}^{2}},\qquad d\ll \mu _{%
\mathrm{d}}^{2}  \label{b3.12}
\end{equation}%
Hence, $\vartheta _{\mathrm{dM}}\approx \sqrt{2d}/\mu _{\mathrm{d}}$. It
means, that the adjacent link is located on the cone of angle $\sqrt{2d}/\mu
_{\mathrm{d}}$, and the whole line $\mathcal{T}_{\mathrm{br}}\left( \sigma _{%
\mathrm{d}}\right) $ has a random shape, because any link wobbles with the
characteristic angle $\sqrt{2d}/\mu _{\mathrm{d}}$. The wobble angle depends
on the space-time distortion $d$ and on the particle mass $\mu _{\mathrm{d}}$%
. The wobble angle is small for the large mass of a particle. The random
displacement of the segment end is of the order $\mu _{\mathrm{d}}\vartheta
_{\mathrm{dM}}=\sqrt{2d}$, i.e. of the same order as the segment width. It
is reasonable, because these two phenomena have the common source: the
space-time distortion $D$.

One should note that the space-time geometry influences the stochasticity of
particle motion illocally in the sense, that the form of the world function (%
\ref{b3.3}) for values of $\sigma _{\mathrm{M}}<\frac{1}{2}\mu _{\mathrm{d}%
}^{2}$ is unessential for the motion stochasticity of the particle of the
mass $\mu _{\mathrm{d}}$.

Such a situation, when the world line of a free particle is stochastic in
the deterministic geometry, and this stochasticity depends on the particle
mass, seems to be rather exotic and incredible. But experiments show that
the motion of real particles of small mass is stochastic indeed, and this
stochasticity increases, when the particle mass decreases. From physical
viewpoint a theoretical foundation of the stochasticity is desirable, and
some researchers invent stochastic geometries, noncommutative geometries and
other exotic geometrical constructions, to obtain the quantum stochasticity.
But in the Riemannian space-time geometry the particle motion does not
depend on the particle mass, and in the framework of the Riemannian
space-time geometry it is difficult to explain the quantum stochasticity by
the space-time geometry properties. The distorted geometry $\mathcal{G}_{%
\mathrm{d}}$ explains freely the stochasticity and its dependence on the
particle mass. Besides, at proper choice of the distortion $d$ the
statistical description of stochastic $\mathcal{T}_{\mathrm{br}}$ leads to
the quantum description (in terms of the Schr\"{o}dinger equation) \cite{R91}%
. To do this, it is sufficient to set 
\begin{equation}
d=\frac{\hbar }{2bc}  \label{b3.12a}
\end{equation}
where $\hbar $ is the quantum constant, $c$ is the speed of the light, and $%
b $ is some universal constant, connecting the geometrical mass $\mu $ with
the usual particle mass $m$ by means of the relation $m=b\mu $. In other
words, the distorted space-time geometry (\ref{b3.3}) is closer to the real
space-time geometry, than the Minkowski geometry $\mathcal{G}_{\mathrm{M}}$.

Further development of the statistical description of geometrical
stochasticity leads to a creation of the model conception of quantum
phenomena (MCQP), which relates to the conventional quantum theory
approximately in the same way as the statistical physics relates to the
axiomatic thermodynamics. MCQP is the well defined relativistic conception
with effective methods of investigation \cite{R03}, whereas the conventional
quantum theory is not well defined, because it uses incorrect space-time
geometry, whose incorrectness is compensated by additional hypotheses
(quantum principles). Besides, it has problems with application of the
nonrelativistic quantum mechanical technique to the description of
relativistic phenomena.

The geometry $\mathcal{G}_{\mathrm{d}}$, as well as the Minkowski geometry
are homogeneous geometries, because the world function $\sigma _{\mathrm{d}}$
is invariant with respect to all coordinate transformations, with respect to
which the world function $\sigma _{\mathrm{M}}$ is invariant.

\end{document}